\documentclass[11pt]{article}
\usepackage{amsfonts}
\usepackage{mathrsfs}
\usepackage{amscd}
\usepackage{amsmath,amsfonts,amssymb,amscd}
\usepackage{indentfirst,graphics,epsfig,psfrag}
\input{epsf}
\usepackage{ifpdf}
\usepackage{enumerate}
\usepackage{appendix}
\usepackage{enumerate}
\usepackage{lineno}

\makeatletter \@addtoreset{figure}{section} \makeatother
\makeatletter
\long\def\@makecaption#1#2{%
   \vskip 10\p@
   \setbox\@tempboxa\hbox{{#1}\ \ #2}%
   \ifdim \wd\@tempboxa >\hsize

       {#1}\ \ #2\par
   \else
       \hbox to\hsize{\hfil\box\@tempboxa\hfil}%
   \fi}
\makeatother

\newtheorem{thm}{Theorem}[section]
\newtheorem{cor}{Corollary}[section]
\newtheorem{lem}{Lemma}[section]

\newtheorem{obs}{Observation}[section]
\newtheorem{pro}{Proposition}[section]

\newcommand{\qed}{{\hfill\rule{3pt}{7pt}}}

\def\qed{\hfill \rule{4pt}{7pt}}

\begin{document}
\title{\textbf{Pendant-tree connectivity
of line graphs}\footnote{Supported by the National Science
Foundation of China (Nos. 11551001, 11161037, 11461054) and the
Science Found of Qinghai Province (No. 2014-ZJ-907).}}
\author{
\small Yaping Mao$^{1,2}$\footnote{E-mail: maoyaping@ymail.com}\\[0.2cm]
\small $^{1}$Department of Mathematics, Qinghai Normal University,\\[0.1mm]
\small $^{2}$Key Laboratory of IOT of Qinghai Province,\\[0.1mm]
\small Xining, Qinghai 810008, China}
\date{}
\maketitle
\begin{abstract}
The concept of pendant-tree connectivity, introduced by Hager in
1985, is a generalization of classical vertex-connectivity. In this
paper, we study pendant-tree connectivity of line graphs.\\[2mm]
\end{abstract}

\section{Introduction}

A processor network is expressed as a graph, where a node is a
processor and an edge is a communication link. Broadcasting is the
process of sending a message from the source node to all other nodes
in a network. It can be accomplished by message dissemination in
such a way that each node repeatedly receives and forwards messages.
Some of the nodes and/or links may be faulty. However, multiple
copies of messages can be disseminated through disjoint paths. We
say that the broadcasting succeeds if all the healthy nodes in the
network finally obtain the correct message from the source node
within a certain limit of time. A lot of attention has been devoted
to fault-tolerant broadcasting in networks \cite{Fragopoulou,
Hedetniemi, Jalote, Ramanathan}. In order to measure the ability of
fault-tolerance, the above path structure connecting two nodes are
generalized into some tree structures connecting more than two
nodes, see \cite{Ku, LLSun, LM1}. To show these generalizations
clearly, we must state from the connectivity in graph theory. We
divide our introduction into the following four subsections to state
the motivations and our results of this paper.

\subsection{Connectivity and $k$-connectivity}

All graphs considered in this paper are undirected, finite and
simple. We refer to the book \cite{bondy} for graph theoretical
notation and terminology not described here. For a graph $G$, let
$V(G)$, $E(G)$, $e(G)$, $L(G)$ and $\delta(G)$ denote the set of
vertices, the set of edges, the size, the line graph and the minimum
degree of $G$, respectively. In the sequel, let $K_{s,t}$, $K_{n}$
and $P_n$ denote the complete bipartite graph of order $s+t$ with
part sizes $s$ and $t$, complete graph of order $n$, and path of
order $n$, respectively. For any subset $X$ of $V(G)$, let $G[X]$
denote the subgraph induced by $X$, and $E[X]$ the edge set of
$G[X]$. For two subsets $X$ and $Y$ of $V(G)$ we denote by
$E_G[X,Y]$ the set of edges of $G$ with one end in $X$ and the other
end in $Y$. If $X=\{x\}$, we simply write $E_G[x,Y]$ for
$E_G[\{x\},Y]$.

Connectivity is one of the most basic concepts of graph-theoretic
subjects, both in combinatorial sense and the algorithmic sense. It
is well-known that the classical connectivity has two equivalent
definitions. The \emph{connectivity} of $G$, written $\kappa(G)$, is
the minimum order of a vertex set $S\subseteq V(G)$ such that $G-S$
is disconnected or has only one vertex. We call this definition the
`cut' version definition of connectivity. A well-known theorem of
Whitney \cite{Whitney} provides an equivalent definition of
connectivity, which can be called the `path' version definition of
connectivity. For any two distinct vertices $x$ and $y$ in $G$, the
\emph{local connectivity} $\kappa_{G}(x,y)$ is the maximum number of
internally disjoint paths connecting $x$ and $y$. Then
$\kappa(G)=\min\{\kappa_{G}(x,y)\,|\,x,y\in V(G),x\neq y\}$ is
defined to be the \emph{connectivity} of $G$. For connectivity,
Oellermann gave a survey paper on this subject; see
\cite{Oellermann2}.

Although there are many elegant and powerful results on connectivity
in graph theory, the basic notation of classical connectivity may
not be general enough to capture some computational settings. So
people want to generalize this concept. For the `cut' version
definition of connectivity, we find the above minimum vertex set
without regard the number of components of $G-S$. Two graphs with
the same connectivity may have differing degrees of vulnerability in
the sense that the deletion of a vertex cut-set of minimum
cardinality from one graph may produce a graph with considerably
more components than in the case of the other graph. For example,
the star $K_{1,n}$ and the path $P_{n+1}\ (n\geq 3)$ are both trees
of order $n+1$ and therefore connectivity $1$, but the deletion of a
cut-vertex from $K_{1,n}$ produces a graph with $n$ components while
the deletion of a cut-vertex from $P_{n+1}$ produces only two
components. Chartrand et al. \cite{Chartrand1} generalized the `cut'
version definition of connectivity. For an integer $k \ (k\geq 2)$
and a graph $G$ of order $n \ (n\geq k)$, the
\emph{$k$-connectivity} $\kappa'_k(G)$ is the smallest number of
vertices whose removal from $G$ of order $n \ (n\geq k)$ produces a
graph with at least $k$ components or a graph with fewer than $k$
vertices. Thus, for $k=2$, $\kappa'_2(G)=\kappa(G)$. For more
details about $k$-connectivity, we refer to \cite{Chartrand1, Day,
Oellermann2, Oellermann3}.

\subsection{Generalized (edge-)connectivity}

The generalized connectivity of a graph $G$, introduced by Hager
\cite{Hager}, is a natural generalization of the `path' version
definition of connectivity. For a graph $G=(V,E)$ and a set
$S\subseteq V(G)$ of at least two vertices, \emph{an $S$-Steiner
tree} or \emph{a Steiner tree connecting $S$} (or simply, \emph{an
$S$-tree}) is a such subgraph $T=(V',E')$ of $G$ that is a tree with
$S\subseteq V'$. Note that when $|S|=2$ an $S$-Steiner tree is just
a path connecting the two vertices of $S$. Two $S$-Steiner trees $T$
and $T'$ are said to be \emph{internally disjoint} if $E(T)\cap
E(T')=\varnothing$ and $V(T)\cap V(T')=S$. For $S\subseteq V(G)$ and
$|S|\geq 2$, the \emph{generalized local connectivity} $\kappa_G(S)$
is the maximum number of internally disjoint $S$-Steiner trees in
$G$, that is, we search for the maximum cardinality of edge-disjoint
trees which include $S$ and are vertex disjoint with the exception
of $S$. For an integer $k$ with $2\leq k\leq n$, \emph{generalized
$k$-connectivity} (or \emph{$k$-tree-connectivity}) is defined as
$\kappa_k(G)=\min\{\kappa_G(S)\,|\,S\subseteq V(G),|S|=k\}$, that
is, $\kappa_k(G)$ is the minimum value of $\kappa_G(S)$ when $S$
runs over all $k$-subsets of $V(G)$. Clearly, when $|S|=2$,
$\kappa_2(G)$ is nothing new but the connectivity $\kappa(G)$ of
$G$, that is, $\kappa_2(G)=\kappa(G)$, which is the reason why one
addresses $\kappa_k(G)$ as the generalized connectivity of $G$. By
convention, for a connected graph $G$ with less than $k$ vertices,
we set $\kappa_k(G)=1$. Set $\kappa_k(G)=0$ when $G$ is
disconnected. Note that the generalized $k$-connectivity and
$k$-connectivity of a graph are indeed different. Take for example,
the graph $H_1$ obtained from a triangle with vertex set
$\{v_1,v_2,v_3\}$ by adding three new vertices $u_1,u_2,u_3$ and
joining $v_i$ to $u_i$ by an edge for $1 \leq i\leq 3$. Then
$\kappa_3(H_1)=1$ but $\kappa'_3(H_1)=2$. There are many results on
the generalized connectivity, see \cite{Chartrand2, LLSun, LL, LLZ,
LM1, LM2, LM3, LM4, LMS, Okamoto}.

As a natural counterpart of the generalized connectivity, we
introduced the concept of generalized edge-connectivity in
\cite{LMS}. For $S\subseteq V(G)$ and $|S|\geq 2$, the {\it
generalized local edge-connectivity} $\lambda(S)$ is the maximum
number of edge-disjoint Steiner trees connecting $S$ in $G$. For an
integer $k$ with $2\leq k\leq n$, the {\it generalized
$k$-edge-connectivity} $\lambda_k(G)$ of $G$ is then defined as
$\lambda_k(G)=\min\{\lambda(S)\,|\,S\subseteq V(G) \ and \ |S|=k\}$.
It is also clear that when $|S|=2$, $\lambda_2(G)$ is just the
standard edge-connectivity $\lambda(G)$ of $G$, that is,
$\lambda_2(G)=\lambda(G)$, which is the reason why we address
$\lambda_k(G)$ as the generalized edge-connectivity of $G$. Also set
$\lambda_k(G)=0$ when $G$ is disconnected. Results on the
generalized edge-connectivity can be found in \cite{LM1, LM4, LMS}.

\subsection{Pendant-tree (edge-)connectivity}

The concept of pendant-tree connectivity \cite{Hager} was introduced
by Hager in 1985, which is specialization of generalized
connectivity (or \emph{$k$-tree-connectivity}) but a generalization
of classical connectivity. For an $S$-Steiner tree, if the degree of
each vertex in $S$ is equal to one, then this tree is called a
\emph{pendant $S$-Steiner tree}. Two pendant $S$-Steiner trees $T$
and $T'$ are said to be \emph{internally disjoint} if $E(T)\cap
E(T')=\varnothing$ and $V(T)\cap V(T')=S$. For $S\subseteq V(G)$ and
$|S|\geq 2$, the \emph{pendant-tree local connectivity} $\tau_G(S)$
is the maximum number of internally disjoint pendant $S$-Steiner
trees in $G$. For an integer $k$ with $2\leq k\leq n$,
\emph{pendant-tree $k$-connectivity} is defined as
$\tau_k(G)=\min\{\tau_G(S)\,|\,S\subseteq V(G),|S|=k\}$. When $k=2$,
$\tau_2(G)=\tau(G)$ is just the connectivity of a graph $G$. For
more details on pendant-tree connectivity, we refer to \cite{Hager,
MaoL}. Clearly, we have
$$
\left\{
\begin{array}{ll}
\tau_k(G)=\kappa_k(G),&\mbox {\rm for}~k=1,2;\\
\tau_k(G)\leq \kappa_k(G),&\mbox {\rm for}~k\geq 3.
\end{array}
\right.
$$

The relation between pendant-tree connectivity and generalized
connectivity are shown in the following Table 2.
\begin{center}
\begin{tabular}{|c|c|c|}
\hline &Pendant tree-connectivity & Generalized connectivity\\[0.1cm]
\cline{1-3}
Vertex subset & $S\subseteq V(G) \ (|S|\geq 2)$ & $S\subseteq V(G) \ (|S|\geq 2)$\\[0.1cm]
\cline{1-3} Set of Steiner trees & $\left\{
\begin{array}{ll}
\mathscr{T}_{S}=\{T_1,T_2,\cdots,T_{\ell}\}\\
S\subseteq V(T_i),\\
d_{T_i}(v)=1~for~every~v\in S\\
E(T_i)\cap E(T_j)=\varnothing,\\
\end{array}
\right.$ & $\left\{
\begin{array}{ll}
\mathscr{T}_{S}=\{T_1,T_2,\cdots,T_{\ell}\}\\
S\subseteq V(T_i),\\
E(T_i)\cap E(T_j)=\varnothing,\\
\end{array}
\right.$\\[0.05cm]
\cline{1-3}
Local parameter & $\tau(S)=\max|\mathscr{T}_{S}|$ & $\kappa(S)=\max|\mathscr{T}_{S}|$\\
\cline{1-3} Global parameter & $\tau_k(G)=\underset{S\subseteq V(G),
|S|=k}{\min} \tau(S)$ & $\kappa_k(G)=\underset{S\subseteq V(G),
|S|=k}{\min}
\kappa(S)$\\[0.05cm]
\cline{1-3}
\end{tabular}
\end{center}
\begin{center}
{Table 2. Two kinds of tree-connectivities}
\end{center}

The following two observations are easily seen.
\begin{obs}\label{obs1-1}
If $G$ is a connected graph, then $\tau_k(G)\leq \mu_k(G)\leq
\delta(G)$.
\end{obs}
\begin{obs}\label{obs1-2}
If $H$ is a spanning subgraph of $G$, then $\tau_k(H)\leq
\tau_k(G)$.
\end{obs}

In {\upshape\cite{Hager}}, Hager derived the following results.
\begin{lem}{\upshape\cite{Hager}}\label{lem1-1}
Let $G$ be a graph. If $\tau_k(G)\geq \ell$, then $\delta(G)\geq
k+\ell-1$.
\end{lem}
\begin{lem}{\upshape\cite{Hager}}\label{lem1-2}
Let $G$ be a graph. If $\tau_k(G)\geq \ell$, then $\kappa(G)\geq
k+\ell-2$.
\end{lem}
\begin{lem}{\upshape\cite{Hager}}\label{lem1-3}
Let $k,n$ be two integers with $3\leq k\leq n$, and let $K_n$ be a
complete graph of order $n$. Then
$$
\tau_k(K_n)=n-k.
$$
\end{lem}

As a natural counterpart of the pendant-tree $k$-connectivity, we
introduced the concept of pendant-tree $k$-edge-connectivity. For
$S\subseteq V(G)$ and $|S|\geq 2$, the {\it pendant-tree local
edge-connectivity} $\mu(S)$ is the maximum number of edge-disjoint
pendant $S$-Steiner trees in $G$. For an integer $k$ with $2\leq
k\leq n$, the {\it pendant-tree $k$-edge-connectivity} $\mu_k(G)$ of
$G$ is then defined as $\mu_k(G)=\min\{\mu(S)\,|\,S\subseteq V(G) \
and \ |S|=k\}$. It is also clear that when $|S|=2$, $\mu_2(G)$ is
just the standard edge-connectivity $\lambda(G)$ of $G$, that is,
$\mu_2(G)=\lambda(G)$.

\subsection{Application background and our results}

In addition to being a natural combinatorial measure, both the
pendant-tree connectivity and the generalized connectivity can be
motivated by its interesting interpretation in practice. For
example, suppose that $G$ represents a network. If one considers to
connect a pair of vertices of $G$, then a path is used to connect
them. However, if one wants to connect a set $S$ of vertices of $G$
with $|S|\geq 3$, then a tree has to be used to connect them. This
kind of tree with minimum order for connecting a set of vertices is
usually called a Steiner tree, and popularly used in the physical
design of VLSI (see \cite{Grotschel1, Grotschel2, Sherwani}) and
computer communication networks (see \cite{Du}). Usually, one wants
to consider how tough a network can be, for the connection of a set
of vertices. Then, the number of totally independent ways to connect
them is a measure for this purpose. The generalized $k$-connectivity
can serve for measuring the capability of a network $G$ to connect
any $k$ vertices in $G$.

Chartrand and Stewart \cite{Steeart} investigated the relation
between the connectivity and edge-connectivity of a graph and its
line graph.
\begin{thm}{\upshape \cite{Steeart}}\label{th1-1}
If $G$ is a connected graph, then

$(1)$ $\kappa(L(G))\geq \lambda(G)$ if $\lambda(G)\geq 2$.

$(2)$ $\lambda(L(G))\geq 2\lambda(G)-2$.

$(3)$ $\kappa(L(L(G)))\geq 2\kappa(G)-2$.
\end{thm}

In Section 2, we investigate the relation between the pendant-tree
$3$-connectivity and pendant-tree $3$-edge-connectivity of a graph
and its line graph.

In their book, Capobianco and Molluzzo \cite{Capobianco}, using
$K_{1,n}$ as their example, note that the difference between the
connectivity of a graph and its line graph can be arbitrarily large.
They then proposed an open problem: Whether for any two integers
$p,q \ (1<p<q)$, there exists a graph $G$ such that $\kappa(G)=p$
and $\kappa(L(G))=q$. In \cite{Bauer}, Bauer and Tindell gave a
positive answer of this problem, that is, for every pair of integers
$p,q \ (1<p<q)$ there is a graph of connectivity $p$ whose line
graph has connectivity $q$.

Note that the difference between the pendant-tree $k$-connectivity
of a graph $G$ and its line graph $L(G)$ can be also arbitrarily
large. Let $n,k$ be two integers with $2\leq k\leq n$, and
$G=K_{1,n}$. Then $L(G)=K_n$, $\tau_k(G)=0$ and $\tau_k(L(G))=n-k$.
In fact, we can consider a similar problem: Whether for any two
integers $p,q$, $1<p<q$, there exists a graph $G$ such that
$\tau_k(G)=p$ and $\tau_k(L(G))=q$. It seem to be not easy to solve
this problem for a general $k$. In this paper, we focus our
attention on the case $k=3$, and give a positive answer of this
problem.

\section{Preliminary}

In \cite{Hager}, Hager showed that $\tau_k$ is monotonically
decreasing for $k$.
\begin{lem}{\upshape \cite{Hager}}\label{lem2-1}
Let $G$ be a graph, and let $k$ be an integer with $k\geq 2$. Then
$$
\tau_k(G)\geq \tau_{k+1}(G).
$$
\end{lem}

By a result in \cite{Hager}, Mao and Lai obtained the following
bounds of $\tau_k(G)$.
\begin{lem}{\upshape \cite{MaoL}}\label{lem2-2}
Let $G$ be a graph of order $n$, and let $k$ be an integer with
$3\leq k\leq n$. Then
$$
\frac{1}{k+1}\log_{2}\kappa(G)\leq \tau_k(G)\leq \kappa(G).
$$
Moreover, the bounds are sharp.
\end{lem}
\begin{pro}{\upshape \cite{MaoLai}}\label{pro2-1}
Let $k,n$ be two integers with $3\leq k\leq n$, and let $G$ be a
graph. Then
$$
0\leq \tau_k(G)\leq n-k.
$$
Moreover, the bounds are sharp.
\end{pro}

For $k=n,n-1,n-2$, the following corollaries are immediate.
\begin{lem}{\upshape \cite{MaoLai}}\label{lem2-3}
Let $G$ be a graph of order $n$. Then $\tau_n(G)=0$ if and only if
$G$ is a graph of order $n$.
\end{lem}
\begin{lem}{\upshape \cite{MaoLai}}\label{lem2-4}
Let $G$ be a connected graph of order $n$. Then

$(1)$ $\tau_{n-1}(G)=1$ if and only if $G$ is a complete graph of
order $n$.

$(2)$ $\tau_{n-1}(G)=0$ if and only if $G$ is not a complete graph
of order $n$.
\end{lem}
\begin{lem}{\upshape \cite{MaoLai}}\label{lem2-5}
Let $G$ be a connected graph of order $n$. Then

$(1)$ $\tau_{n-2}(G)=2$ if and only if $G$ is a complete graph of
order $n$.

$(2)$ $\tau_{n-2}(G)=1$ if and only if $G=K_n\setminus M$ and $1\leq
|M|\leq 2$, where $M$ is a matching of $K_n$ for $n\geq 7$.

$(3)$ $\tau_{n-2}(G)=0$ if and only if $G$ is one of the other
graphs.
\end{lem}

The following results for line graphs can be found in \cite{West}.
\begin{lem}{\upshape \cite{Chang, Hoffman}}\label{lem2-6}
For $n\neq 8$, $L(K_n)$ is the only $(2n-4)$-regular simple graph of
order ${n\choose {2}}$ in which nonadjacent vertices have four
common neighbors and adjacent vertices have $n-2$ common neighbors.
\end{lem}

\begin{lem}{\upshape (\cite{West}, p-283)}\label{lem2-7}
Let $G$ be a $k$-edge-connected simple graph. Then $L(G)$ is
$k$-connected and $(2k-2)$-edge-connected.
\end{lem}

Let $S$ be a set of $k$ vertices of a connected graph $G$, and let
$\mathcal{T}$ be a set of edge-disjoint pendant $S$-Steiner trees.
Then the following observation is immediate.
\begin{obs}\label{obs2-1}
Let $k,n$ be two integers with $3\leq k\leq n$. Let $G$ be a graph
of order $n$, and let $S\subseteq V(G)$ with $|S|=k$. For each $T\in
\mathcal{T}$,
$$
|E(T)\cap E_G[S,\bar{S}]|\geq k,
$$
where $\bar{S}=V(G)-S$.
\end{obs}

By the above result, we can easily derive an upper bound for
pendant-tree $k$-edge-connectivity.
\begin{obs}\label{obs2-2}
For any graph $G$ with order at least $k$,
$$
\mu_k(G)\leq \min_{S\subseteq V(G),
|S|=k}\Big\lfloor\frac{1}{k}|E_G[S,\bar{S}]|\Big\rfloor,
$$
where $S\subseteq V(G)$ with $|S|=k$, and $\bar{S}=V(G)-S$.
Moreover, the bound is sharp.
\end{obs}

\begin{pro}\label{pro2-3}
Let $k,n$ be two integers with $3\leq k\leq n$, and let $K_n$ be a
complete graph of order $n$. Then
$$
\mu_k(K_n)=n-k.
$$
\end{pro}
\begin{pf}
From Observation \ref{obs2-2}, we have $\mu_{k}(K_n)\leq
\min_{S\subseteq V(K_n),
|S|=k}\lfloor\frac{1}{k}|E_{K_n}[S,\bar{S}]|\rfloor=\frac{1}{k}(n-k)k=n-k$.
From Lemma \ref{lem1-3} and Observation \ref{obs1-1}, we have
$\mu_{k}(K_n)\geq \tau_{k}(K_n)=n-k$. So $\mu_k(K_n)=n-k$, as
desired.\qed
\end{pf}

Note that each graph is a spanning subgraph of a complete graph. So
the following result is immediate.
\begin{pro}\label{pro2-4}
Let $k,n$ be two integers with $3\leq k\leq n$, and let $G$ be a
graph. Then
$$
0\leq \mu_k(G)\leq n-k.
$$
\end{pro}

Graphs with $\mu_k(G)=n-k$ are characterized as follows.
\begin{pro}\label{pro2-5}
Let $G$ be a graph of order $n$. Then $\mu_{k}(G)=n-k$ if and only
if $G$ is a complete graph.
\end{pro}
\begin{pf}
If $G$ is a complete graph, then it follows from Proposition
\ref{pro2-3} that $\mu_{k}(G)\geq \tau_{k}(G)\geq n-k$. From
Proposition \ref{pro2-4}, we have $\mu_{k}(G)=n-k$. Conversely, we
suppose $\mu_{k}(G)=n-k$. We claim that $G$ is a complete graph.
Assume, to the contrary, that $G$ is not a complete graph. From
Observation \ref{obs2-2}, we have
$$\mu_{k}(G)\leq \min_{S\subseteq
V(G),
|S|=n-2}\left\lfloor\frac{1}{k}|E_{K_n}[S,\bar{S}]|\right\rfloor\leq
\frac{1}{k}[k(n-k)-1]=n-k-\frac{1}{k},
$$
and hence $\mu_{k}(G)\leq n-k-1$, a contradiction. So $G$ is a
complete graph of order $n$, as desired.
\end{pf}

The following two corollaries are immediate from Propositions
\ref{pro2-4} and \ref{pro2-5}.
\begin{cor}\label{cor2-1}
Let $G$ be a graph of order $n$. Then $\mu_n(G)=0$ if and only if
$G$ is a graph of order $n$.
\end{cor}
\begin{cor}\label{cor2-2}
Let $G$ be a graph of order $n$. Then

$(1)$ $\mu_{n-1}(G)=1$ if and only if $G$ is a complete graph.

$(2)$ $\mu_{n-1}(G)=0$ if and only if $G$ is not a complete graph.
\end{cor}

From the above corollary, we can get the relation between
$\mu_{n-1}(G)$ and $\lambda(G)$.
\begin{pro}\label{pro2-6}
Let $G$ be a graph of order $n$. Then
$$
\mu_{n-1}(G)=\left \lfloor \frac{\lambda(G)}{n-1}\right \rfloor.
$$
\end{pro}
\begin{pf}
If $G$ is a complete graph, then it follows from Corollary
\ref{cor2-2} that $\lambda(G)=n-1$ and $\mu_{n-1}(G)=1$, and hence
$\mu_{n-1}(G)=1=\lfloor \frac{\lambda(G)}{n-1}\rfloor$. If $G$ is
not a complete graph, then it follows from Corollary \ref{cor2-2}
that $0\leq \lambda(G)\leq n-2$ and $\mu_{n-1}(G)=0$, and hence
$\mu_{n-1}(G)=0=\lfloor \frac{\lambda(G)}{n-1}\rfloor$, as
desired.\qed
\end{pf}

From Proposition \ref{pro2-4}, we have $0\leq \mu_{n-2}(G)\leq 2$.
In the following, graphs with $\mu_{n-2}(G)=\ell \ (0\leq \ell\leq
2)$ are characterized.

\begin{lem}\label{lem2-8}
Let $G$ be a connected graph of order $n$. Then $\mu_{n-2}(G)=1$ if
and only if $\bar{G}=K_{1,r}\cup (n-r-1)K_1$, where $1\leq r\leq
n-2$.
\end{lem}
\begin{pf}
Suppose that $\bar{G}=K_{1,r}\cup (n-r-1)K_1$ and $1\leq r\leq n-2$.
From Proposition \ref{pro2-5}, we have $\mu_{n-2}(G)\leq 1$. Let $u$
be the center of $K_{1,r}$. For any $S\subseteq V(G)$ and $|S|=n-2$,
we have $u\in S$ or $u\in \bar{S}$, where $\bar{S}=V(G)-S$. If $u\in
\bar{S}$, then we let $\bar{S}=\{u,v\}$. Since $\bar{G}=K_{1,r}\cup
(n-r-1)K_1$, it follows that the tree induced by the edges in
$E_G[v,S]$ is a pendant $S$-Steiner tree, and hence $\mu(S)\geq 1$.
From now on, we suppose $u\in S$. Set $\bar{S}=\{v,w\}$. Suppose
that $u\in S$. Suppose $uv\notin E(G)$ and $uw\notin E(G)$. Since
$\bar{G}=K_{1,r}\cup (n-r-1)K_1$ where $1\leq r\leq n-2$, it follows
that there exists a vertex in $S$, say $x$, such that $ux\in E(G)$.
Then the tree induced by the edges in $\{ux\}\cup E_G[w,S]-\{uw\}$
is a pendant $S$-Steiner tree, and hence $\mu(S)\geq 1$. Suppose
$uv\notin E(G)$ or $uw\notin E(G)$. Without loss of generality, let
$uv\notin E(G)$. Then $uw\in E(G)$, and the tree induced by the
edges in $E_G[w,S]$ is a pendant $S$-Steiner tree, and hence
$\mu(S)\geq 1$. From the argument, we conclude that $\mu(S)\geq 1$
for any $S\subseteq V(G)$ and $|S|=n-2$. From the arbitrariness of
$S$, we have $\mu_{n-2}(G)\geq 1$, and hence $\mu_{n-2}(G)=1$.

Conversely, we suppose $\mu_{n-2}(G)=1$. Then we have the following
claim.

\noindent \textbf{Claim 1.} $\bar{G}$ does not contain a $P_4$ as
its subgraph.

\noindent \textbf{Proof of Claim 1:} Assume, to the contrary, that
$\bar{G}$ contains a $P_4=u_1u_2u_3u_4$ as its subgraph. Choose
$S=V(G)-\{u_2u_3\}$. Since $u_iu_{i+1}\notin E(G) \ (1\leq i\leq
3)$, it follows that there is no $S$-Steiner tree in $G$, a
contradiction.

From Claim 1, $\bar{G}$ does not contain cycles of order at least
$4$ as its subgraph. Furthermore, we have the following claim.

\noindent \textbf{Claim 2.} $\bar{G}$ does not contain cycles.

\noindent \textbf{Proof of Claim 2:} Assume, to the contrary, that
$\bar{G}$ contains a cycle. From Claim 1, this cycle is a triangle,
say $C_3=uvw$. Choose $S=V(G)-\{u,w\}$. Since $uv,uw,vw\notin E(G)$,
it follows that there is no pendant $S$-Steiner tree in $G$, a
contradiction.

From Claim 2, $\bar{G}$ is a tree. From Claim 1, $\bar{G}$ is a
star. Set $\bar{G}=K_{1,r}$. Since $\mu_{n-2}(G)=1$, it follows that
$G$ is connected, and hence $1\leq r\leq n-2$, as desired. \qed
\end{pf}

From Lemma \ref{lem2-8} and Proposition \ref{pro2-5}, the following
result is easily seen.
\begin{pro}\label{pro2-7}
Let $G$ be a connected graph of order $n$. Then

$(1)$ $\mu_{n-2}(G)=2$ if and only if $G$ is a complete graph of
order $n$.

$(2)$ $\mu_{n-2}(G)=1$ if and only if $\bar{G}=K_{1,r}\cup
(n-r-1)K_1$ where $1\leq r\leq n-2$.

$(3)$ $\mu_{n-2}(G)=0$ if and only if $e(\bar{G})\geq 1$ and
$\bar{G}\neq K_{1,r}\cup (n-r-1)K_1$ and $1\leq r\leq n-2$.
\end{pro}

From Proposition \ref{pro2-7}, we can set up the relation between
$\mu_{n-2}(G)$ and $\lambda(G)$.
\begin{pro}\label{pro2-8}
Let $G$ be a graph of order $n$. Then
$$
\mu_{n-2}(G)=\left\{
\begin{array}{ll}
\left \lceil \frac{\lambda(G)}{n-2}\right \rceil &\mbox{\rm if}~\lambda(G)=n-1,\\
&{\rm or}~\bar{G}=K_{1,r}\cup
(n-r-1)K_1~{\rm where}~1\leq r\leq n-2,\\[0.2cm]
\left \lfloor \frac{\lambda(G)}{n-2}\right \rfloor &\mbox{\rm
otherwise}.
\end{array}
\right.
$$
\end{pro}
\begin{pf}
If $\lambda(G)=n-1$, then $G$ is a complete graph, and hence
$\mu_{n-2}(G)=2=\lceil \frac{\lambda(G)}{n-2}\rceil$. If
$\bar{G}=K_{1,r}\cup (n-r-1)K_1$ where $1\leq r\leq n-2$, then it
follows from Proposition \ref{pro2-7} that $\mu_{n-2}(G)=1=\lceil
\frac{\lambda(G)}{n-2}\rceil$. For other cases, from Proposition
\ref{pro2-7}, we have
$\mu_{n-2}(G)=0=\lfloor\frac{\lambda(G)}{n-2}\rfloor$.\qed
\end{pf}

\begin{cor}\label{cor2-3}
Let $G$ be a graph of order $n$. Then
$$
\left \lfloor \frac{\lambda(G)}{n-2}\right \rfloor\leq
\mu_{n-2}(G)\leq \left \lceil\frac{\lambda(G)}{n-2}\right \rceil.
$$
\end{cor}

In \cite{MaoL}, Mao and Lai got the following results.
\begin{lem}{\upshape \cite{MaoL}}\label{lem2-9}
Let $G$ be a connected graph with $n$ vertices. Then
$$
\frac{1}{12}\kappa(G)-\frac{1}{2}\leq \tau_3(G)\leq
\frac{2}{3}\kappa(G).
$$
Moreover, the lower bound is sharp.
\end{lem}

\begin{lem}{\upshape \cite{MaoL}}\label{lem2-10}
Let $G$ be a graph of order $n$. Then
$$
\tau_{n-2}(G)=\left\{
\begin{array}{ll}
\left \lceil \frac{\kappa(G)}{n-2}\right \rceil &\mbox{\rm if}~\kappa(G)=n-1,\\
&{\rm or}~\kappa(G)=n-2~{\rm and}~\bar{G}=iK_2\cup (n-2i)K_1 \ (i=1,2),\\[0.2cm]
\left \lfloor \frac{\kappa(G)}{n-2}\right \rfloor &\mbox{\rm
otherwise}.
\end{array}
\right.
$$
\end{lem}

\begin{lem}\label{lem2-11}
Let $G$ be a graph of order $n$. Then
$$
\tau_{n-1}(G)=\left \lfloor \frac{\kappa(G)}{n-1}\right \rfloor.
$$
\end{lem}

\section{Pendant-tree connectivity of
line graphs}

For pendant-tree $k$-connectivity, we have the following:
\begin{pro}\label{pro3-1}
If $G$ is a connected graph, then

$(1)$ $\frac{1}{k+1}\log_{2}\mu_k(G)\leq \tau_k(L(G))$.

$(2)$ $\mu_k(L(G))\geq \frac{1}{k+1}\log_{2}\mu_k(G)$.

$(3)$ $\tau_k(L(L(G))\geq
\frac{1}{k+1}\left(1+\log_{2}(\tau_k(G)-1)\right)$.
\end{pro}
\begin{pf}
For $(1)$, from Lemma \ref{lem2-2} and Theorem \ref{th1-1}, we have
$$
\tau_k(L(G))\geq \frac{1}{k+1}\log_{2}\kappa(L(G))\geq
\frac{1}{k+1}\log_{2}\lambda(G)\geq \frac{1}{k+1}\log_{2}\mu_k(G).
$$

For $(2)$, from $(1)$ of this proposition, we have
$$
\mu_k(L(G))\geq \tau_k(L(G))\geq \frac{1}{k+1}\log_{2}\mu_k(G).
$$

For $(3)$, from $(1)$ of this proposition and Lemma \ref{lem2-2}, we
have
$$
\tau_k(L(L(G)))\geq \frac{1}{k+1}\log_{2}\kappa(L(L(G)))\geq
\frac{1}{k+1}\log_{2}(2\kappa(G)-2)\geq
\frac{1}{k+1}\left(1+\log_{2}(\tau_k(G)-1)\right).
$$
\end{pf}

As we have seen, the above bounds are relatively rough. In the
following, we try to improve the bounds for $k=n,n-1,n-2,3$.
\begin{pro}\label{pro3-2}
Let $G$ be a connected graph of order $n$. Then

$(1)$ $\mu_n(G)\leq \tau_n(L(G))$.

$(2)$ $\mu_n(L(G))\geq \mu_n(G)$.

$(3)$ $\tau_n(L(L(G))\geq \tau_n(G)$.
\end{pro}
\begin{pf}
For $(1)$, we set $|E(G)|=m$. If $m\geq n$, then it follows from
Corollary \ref{cor2-2} and Lemma \ref{lem2-1} that
$\mu_n(G)=0=\tau_m(L(G))\leq \tau_n(L(G))$. If $m=n-1$, then
$\mu_n(G)=0\leq \tau_n(L(G))$. For $(2)$, from $(1)$, we have
$\mu_n(L(G))\geq \tau_n(L(G))\geq \mu_n(G)$. For $(3)$, from $(1)$
and $(2)$, we have $\tau_n(L(L(G))\geq \mu_n(L(G))\geq \mu_n(G)\geq
\tau_n(G)$. \qed
\end{pf}

\begin{pro}\label{pro3-3}
Let $G$ be a connected graph of order $n \ (n\geq 3)$. If $G$ is
$2$-edge-connected, then

$(1)$ $\left\lceil\frac{2\mu_{n-1}(G)}{n}\right\rceil-1\leq
\tau_{n-1}(L(G))$;

$(2)$ $\mu_{n-1}(L(G))\geq \left\lfloor
\frac{4\mu_{n-1}(G)-4}{(n+1)(n-2)}\right\rfloor$;

$(3)$ $\tau_{n-1}(L(L(G))\geq \left\lfloor
\frac{4\tau_{n-1}(G)-4}{(n+1)(n-2)}\right\rfloor$.
\end{pro}
\begin{pf}
For $(1)$, we set $|E(G)|=m$. Since $G$ is $2$-connected, it follows
that $m\geq n$. From Lemma \ref{lem2-1} and Theorem \ref{th1-1} that
\begin{eqnarray*}
\mu_{n-1}(G)&=&\left\lfloor \frac{\lambda(G)}{n-1}\right\rfloor\leq
\left\lfloor \frac{\kappa(L(G))}{n-1}\right\rfloor\leq
\frac{m}{n-1}\frac{\kappa(L(G))}{m}\\[0.1cm]
&\leq&\frac{m}{n-1}\left(\left\lfloor\frac{\kappa(L(G))}{m}\right\rfloor+1\right)=\frac{m}{n-1}\left(\tau_{m-1}(G)+1\right)\\[0.1cm]
&\leq&\frac{n}{2}\left(\tau_{n-1}(G)+1\right),
\end{eqnarray*}
and hence
$$
\tau_{n-1}(L(G))\geq
\left\lceil\frac{2\mu_{n-1}(G)}{n}\right\rceil-1.
$$

For $(2)$, from Lemma \ref{lem2-1}, Theorem \ref{th1-1} and
Proposition \ref{pro2-5}, we have
$$
\mu_{n-1}(L(G))\geq \mu_{m}(L(G))=\left\lfloor
\frac{\lambda(L(G))}{m-1}\right\rfloor\geq \left\lfloor
\frac{2\lambda(G)-2}{{n\choose 2}-1}\right\rfloor\geq \left\lfloor
\frac{4\mu_{n-1}(G)-4}{(n+1)(n-2)}\right\rfloor.
$$

For $(3)$, from Lemmas \ref{lem2-1}, \ref{lem2-11} and Theorem
\ref{th1-1}, we have
$$
\tau_{n-1}(L(L(G))\geq \tau_{m}(L(L(G))=\left\lfloor
\frac{\kappa(L(L(G))}{m-1}\right\rfloor\geq \left\lfloor
\frac{2\kappa(G)-2}{{n\choose 2}-1}\right\rfloor\geq \left\lfloor
\frac{4\tau_{n-1}(G)-4}{(n+1)(n-2)}\right\rfloor.
$$
\end{pf}

\begin{pro}\label{pro3-3}
Let $G$ be a connected graph of order $n$. If $G$ is
$2$-edge-connected, then

$(1)$ $\left\lfloor \frac{2\mu_{n-2}(G)}{n^2-n-4}\right\rfloor\leq
\tau_{n-2}(L(G))$.

$(2)$ $\mu_{n-2}(L(G))\geq \left\lfloor
\frac{4\mu_{n-2}(G)-4}{n^2-n-4}\right\rfloor$.

$(3)$ $\tau_{n-2}(L(L(G))\geq \left\lfloor
\frac{4\tau_{n-2}(G)-4}{n^2-n-4}\right\rfloor$.
\end{pro}
\begin{pf}
$(1)$ From Theorem \ref{th1-1}, Proposition \ref{pro2-8} and Lemma
\ref{lem2-11}, we have
$$
\tau_{n-2}(L(G))\geq \tau_{m-2}(L(G))\geq \left\lfloor
\frac{\kappa(L(G))}{m-2}\right\rfloor\geq \left\lfloor
\frac{\lambda(G)}{m-2}\right\rfloor\geq \left\lfloor
\frac{2\mu_{n-2}(G)}{n^2-n-4}\right\rfloor.
$$

For $(2)$, from $(1)$ of this proposition and Proposition
\ref{pro2-5}, we have
$$
\mu_{n-1}(L(G))\geq \mu_{m}(L(G))\geq \left\lfloor
\frac{\lambda(L(G))}{m-2}\right\rfloor\geq \left\lfloor
\frac{2\lambda(G)-2}{{n\choose 2}-2}\right\rfloor\geq \left\lfloor
\frac{4\mu_{n-2}(G)-4}{n^2-n-4}\right\rfloor,
$$
where $m$ is the size of $G$.

For $(3)$, from $(1)$ of this proposition and Corollary
\ref{cor2-2}, we have
$$
\tau_{n-2}(L(L(G))\geq \tau_{m}(L(L(G))=\left\lfloor
\frac{\kappa(L(L(G))}{m-2}\right\rfloor\geq \left\lfloor
\frac{2\kappa(G)-2}{{n\choose 2}-2}\right\rfloor\geq \left\lfloor
\frac{4\tau_{n-1}(G)-4}{n^2-n-4}\right\rfloor.
$$
\end{pf}

For pendant-tree $3$-connectivity, we have the following:
\begin{thm}\label{th3-1}
Let $G$ be a connected graph. Then

$(1)$ $\mu_3(G)\leq \tau_3(L(G))$.

$(2)$ $\mu_3(L(G))\geq \frac{1}{12}\mu_3(G)-\frac{1}{2}$.

$(3)$ $\tau_3(L(L(G))\geq \frac{1}{4}\tau_3(G)-\frac{2}{3}$.
\end{thm}
\begin{pf}
For $(1)$, let $e_1,e_2,e_3$ be three arbitrary distinct vertices of
the line graph of $G$ such that $\mu_3(G)=\ell$ with $\ell \geq 1$.
Let $e_1=v_1v_1'$, $e_2=v_2v_2'$ and $e_3=v_3v_3'$ be those edges of
$G$ corresponding to the vertices $e_1,e_2,e_3$ in $L(G)$,
respectively.
\begin{figure}[!hbpt]
\begin{center}
\includegraphics[scale=0.80000]{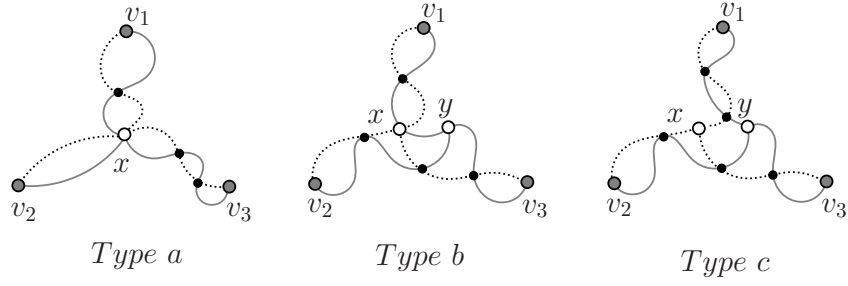}
\end{center}
\caption{Three possible types of $T_i\cup T_j$.}\label{fig3-1}
\end{figure}

Consider three distinct vertices of the six end-vertices of
$e_1,e_2,e_3$. Without loss of generality, let $S=\{v_1,v_2,v_3\}$
be three distinct vertices. Since $\mu_3(G)=\ell$, there exist
$\ell$ edge-disjoint pendant $S$-Steiner trees in $G$, say
$T_1,T_2,\cdots,T_{\ell}$. We define a minimal $S$-Steiner tree $T$
as an $S$-Steiner tree whose subtree obtained by deleting any edge
of $T$ does not connect $S$.

By choosing any two edge-disjoint minimal $S$-Steiner trees $T_i$
and $T_j \ (1\leq i,j\leq \ell)$ in $G$, we will show that the trees
$T_i'$ and $T_j'$ corresponding to $T_i$ and $T_j$ in $L(G)$ are
internally disjoint pendant $S$-Steiner trees. It is easy to see
that $T_i\cup T_j$ has three possible types, as shown in Figure
\ref{fig3-1}. Since $T_i$ and $T_j$ are edge-disjoint in $G$, we can
find internally disjoint pendant Steiner trees $T_i'$ and $T_j'$
connecting $e_1,e_2,e_3$ in $L(G)$. We give an example of Type $a$;
see Figure \ref{fig3-2}. So $\tau_3(L(G))\geq \ell$, as desired.

For $(2)$, from Lemma \ref{lem2-9} and Theorem \ref{th1-1}, we have
$$
\mu_3(L(G))\geq \tau_3(L(G))\geq
\frac{1}{12}\kappa(L(G))-\frac{1}{2}\geq
\frac{1}{12}\lambda(G)-\frac{1}{2}\geq
\frac{1}{12}\mu_3(G)-\frac{1}{2}.
$$

For $(3)$, from Lemma \ref{lem2-9} and Theorem \ref{th1-1}, we have
\begin{eqnarray*}
\tau_3(L(L(G)))&\geq& \frac{1}{12}\kappa(L(L(G)))-\frac{1}{2}\geq \frac{1}{12}(2\kappa(G)-2)-\frac{1}{2}\\
&=&\frac{1}{6}\kappa(G)-\frac{2}{3}\geq
\frac{1}{4}\tau_3(G)-\frac{2}{3}.
\end{eqnarray*}
\begin{figure}[!hbpt]
\begin{center}
\includegraphics[scale=0.80000]{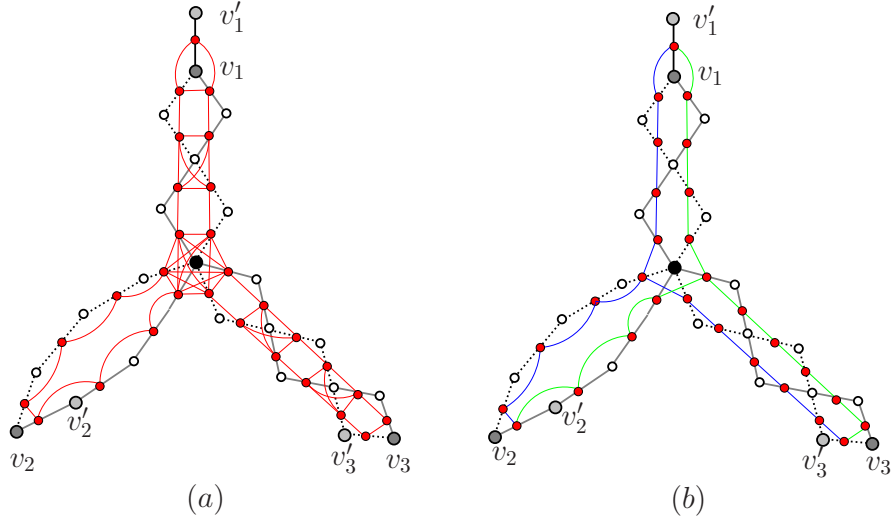}
\end{center}
\caption{(a): An example for $T_i$ and $T_j$ connecting $S$ and
their line graphs; (b): An example for $T_i'$ and $T_j'$
corresponding to $T_i$ and $T_j$.}\label{fig3-2}
\end{figure}
\end{pf}

\section{Graphs with prescribed
pendant-tree connectivity and pendant-tree edge-connectivity}

In {\upshape\cite{Hager}}, Hager obtained the following result.
\begin{lem}{\upshape\cite{Hager}}\label{lem4-1}
Let $K_{a,b}$ be a complete bipartite graph with $a+b$ vertices.
Then
$$
\tau_k(K_{a,b})=\max\{\min\{a-k+1, b-k+1\},0\}.
$$
\end{lem}

The following corollary is immediate from the above lemma.
\begin{cor}\label{cor4-1}
Let $a,b$ be two integers with $2\leq a \leq b$, and $K_{a,b}$
denote a complete bipartite graph with a bipartition of sizes $a$
and $b$, respectively. Then
$$
\tau_3(K_{a,b})=a-2.
$$
\end{cor}

Mao and Lai obtained the following result.
\begin{lem}{\upshape\cite{MaoL}}\label{lem4-2}
Let $G$ be a connected graph with minimum degree $\delta$. Then
$$
\mu_k(G)\leq \delta(G)-k+1.
$$
\end{lem}

In this section, we consider the above problem for the case $p\geq
2q$. Let us put our attention on the complete bipartite graph
$G=K_{p+2,q-p+2}$. Since $q\geq 2p$, it follows that $q-p+2\geq
p+2$. From Corollary \ref{cor4-1},
$\kappa_3(G)=\kappa_3(K_{p+2,q-p+2})=p$. Now we turn to consider the
line graph of the complete bipartite graph $G=K_{p+2,q-p+2}$.

Recall that the Cartesian product (also called the {\em square
product}) of two graphs $G$ and $H$, written as $G\Box H$, is the
graph with vertex set $V(G)\times V(H)$, in which two vertices
$(u,v)$ and $(u',v')$ are adjacent if and only if $u=u'$ and
$(v,v')\in E(H)$, or $v=v'$ and $(u,u')\in E(G)$. Clearly, the
Cartesian product is commutative, that is, $G\Box H\cong H\Box G$.
The following lemma is easily seen.
\begin{lem}{\upshape \cite{West}}\label{lem4-3}
For a complete bipartite graph $K_{r,s}$, $L(K_{r,s})=K_r\Box K_s$.
\end{lem}

From the above lemma, $L(G)=L(K_{p+2,q-p+2})=K_{p+2}\Box K_{q-p+2}$.
In order to obtain the exact value of
$\tau_3(L(G))=\tau_3(K_{p+2}\Box K_{q-p+2})$, we consider to
determine the exact value of the Cartesian product of two complete
graphs.

Before proving Lemma \ref{lem4-4}, we define some notation. Let $G$
and $H$ be two graphs with $V(G)=\{u_1,u_2,\ldots,u_r\}$ and
$V(H)=\{v_1,v_2,\ldots,v_s\}$, respectively. Then $V(G\Box
H)=\{(u_i,v_j)\,|\,1\leq i\leq n, \ 1\leq j\leq m\}$. For $v\in
V(H)$, we use $G(v)$ to denote the subgraph of $G\Box H$ induced by
the vertex set $\{(u_i,v)\,|\,1\leq i\leq n\}$. Similarly, for $u\in
V(G)$, we use $H(u)$ to denote the subgraph of $G\Box H$ induced by
the vertex set $\{(u,v_j)\,|\,1\leq j\leq m\}$.
\begin{lem}\label{lem4-4}
Let $r,s$ be two integers. Then
$$
\tau_3(K_r\Box K_s)=r+s-4.
$$
\end{lem}
\begin{pf}
Let $G=K_r$ and $H=K_s$. Set $V(G)=\{u_1,u_2,\cdots,u_r\}$ and
$V(H)=\{v_1,v_2,\cdots,v_s\}$. From Lemma \ref{lem1-3}, we have
$\tau_3(G)=\tau_3(K_r)=r-3$ and $\tau_3(H)=\tau_3(K_s)=s-3$. On one
hand, we assume $x,y\in V(G(v_1))$. Then $d_{G\Box H}(x)=d_{G\Box
H}(y)=\delta(G\Box H)=r+s-2$. Since $xy\in E(G\Box H)$, it follows
that $\tau_3(G\Box H)\leq \delta(G\Box H)-2\leq r+s-4$.

On the other hand, we will show that $\tau_3(G\Box H)\geq r+s-4$. We
need to show that for any $S=\{x,y,z\}\subseteq V(G\Box H)$, there
exist $r+s-4$ internally disjoint pendant $S$-Steiner trees. We
complete our proof by the following three cases.

\textbf{Case 1}. $x,y,z$ belongs to the same $V(H(u_i)) \ (1\leq
i\leq r)$.

Without loss of generality, we assume $x,y,z\in V(H(u_1))$. Since
$\kappa_3(H)=s-3$, there exist $s-3$ internally disjoint pendant
$S$-Steiner trees $T_1,T_2,\cdots,T_{s-3}$ in $H(u_1)$. Let
$x_j,y_j,z_j$ be the vertices corresponding to $x,y,z$ in $H(u_j) \
(2\leq j\leq r)$. Then the trees $T_{j}'$ induced by the edges in
$\{xx_j,yy_j,zz_j,x_jy_j,y_jz_j\} \ (2\leq j\leq r)$ are $r-1$
internally disjoint pendant $S$-Steiner trees. These trees together
with the trees $T_1,T_2,\cdots,T_{s-3}$ are $r+s-4$ internally
disjoint pendant $S$-Steiner trees.

\textbf{Case 2}. Only two vertices of $\{x,y,z\}$ belong to some
copy $H(u_i)$.

We may assume $x,y\in V(H(u_1))$, $z\in V(H(u_2))$. Let $x',y'$ be
the vertices corresponding to $x,y$ in $H(u_2)$, and let $z'$ be the
vertex corresponding to $z$ in $H(u_1)$. Without loss of generality,
let $V(H(u_1))=\{v_1,v_2,\cdots,v_s\}$ and $V(H(u_2))=\{v_1',v_2',
\cdots,v_s'\}$.

Suppose $z'\not\in \{x,y\}$. Without loss of generality, let
$\{x,y,z'\}=\{v_1,v_2,v_3\}$ in $H(u_1)$. Then the tree $T_1$
induced by the edges in $\{zz',xz',yz'\}$ and $T_{i-2}$ induced by
the edges in $\{xv_i,yv_i,v_iv_i',v_i'z\} \ (4\leq i\leq s)$ are
$s-2$ internally disjoint pendant $S$-Steiner trees. Let
$x_j,y_j,z_j$ be the vertices corresponding to $x,y,z'$ in $H(u_j) \
(3\leq j\leq r)$. The the trees $T_{j}'$ induced by the edges in
$\{xx_j,yy_j, zz_j,x_jz_j,y_jz_j\} \ (3\leq j\leq s)$ are $r-2$
internally disjoint pendant $S$-Steiner trees. These trees together
with the trees $T_1,T_2,\cdots,T_{r-2}$ are $r+s-4$ internally
disjoint pendant $S$-Steiner trees.

Suppose $z'\in \{x,y\}$. Without loss of generality, assume $z'=y$.
Without loss of generality, let $\{x,y\}=\{v_1,v_2\}$ in $H(u_1)$.
Then the trees $T_{i-2}$ induced by the edges in $\{xv_i,
yv_i,v_iv_i',v_i'z\} \ (3\leq i\leq s)$ are $s-2$ internally
disjoint pendant $S$-Steiner trees. Let $x_j,y_j$ be the vertices
corresponding to $x,y$ in $H(u_j) \ (3\leq j\leq r)$. Then the trees
$T_{j}'$ induced by the edges in $\{xx_j,yy_j, zy_j,x_jy_j\} \
(3\leq j\leq r)$ are $r-2$ internally disjoint pendant $S$-Steiner
trees. These trees together with the trees $T_1,T_2,\cdots,T_{r-2}$
are $r+s-4$ internally disjoint pendant $S$-Steiner trees.

\textbf{Case 3}. $x,y,z$ are contained in distinct $H(u_i)$s.

We may assume that $x\in V(H(u_1))$, $y\in V(H(u_2))$, $z\in
V(H(u_3))$. Let $y',z'$ be the vertices corresponding to $y,z$ in
$H(u_1)$, $x',z''$ be the vertices corresponding to $x,z$ in
$H(u_2)$ and $x'',y''$ be the vertices corresponding to $x,y$ in
$H(u_3)$. Without loss of generality, let
$V(H(u_1))=\{v_1,v_2,\cdots,v_s\}$,
$V(H(u_2))=\{v_1',v_2',\cdots,v_s'\}$,
$V(H(u_3))=\{v_1'',v_2'',\cdots,v_s''\}$.

Suppose that $x,y',z'$ are distinct vertices in $H(u_1)$. Without
loss of generality, let $\{x,y',z'\}=\{v_1,v_2,v_3\}$ in $H(u_1)$,
$\{x',y,z''\}=\{v_1',v_2',v_3'\}$ in $H(u_2)$ and
$\{x'',y'',z\}=\{v_1'',v_2'',v_3''\}$ in $H(u_3)$. Then the tree
$T_1$ induced by the edges in $\{xx',x'y,x'x'', x''z\}$, the tree
$T_2$ induced by the edges in $\{xz',z'z'', yz'',z''z\}$, the tree
$T_3$ induced by the edges in $\{xy', yy', yy'',y''z\}$ and the
trees $T_i$ induced by the edges in $\{xu_i,
yu_i',u_i''z,u_iu_i',u_i'u_i''\} \ (4\leq i\leq s)$ are $s$
internally disjoint pendant $S$-Steiner trees. Let $x_j,y_j,z_j$ be
the vertices corresponding to $x,y,z$ in $H(u_j) \ (4\leq j\leq r)$.
The the trees $T_{j}'$ induced by the edges in $\{xx_j,yy_j,
zz_j,x_jy_j,y_jz_j\} \ (4\leq j\leq r)$ are $r-3$ internally
disjoint pendant $S$-Steiner trees. These trees together with the
trees $T_1,T_2,\cdots,T_{r}$ are $r+s-3$ internally disjoint pendant
$S$-Steiner trees.

Suppose that two of $x, y',z'$ are the same vertex in $H(u_1)$.
Without loss of generality, let $y'=z'$, $\{x,y'\}=\{v_1,v_2\}$ in
$H(u_1)$, $\{x',y\}=\{v_1',v_2'\}$ in $H(u_2)$ and
$\{x'',z\}=\{v_1'',v_2''\}$ in $H(u_3)$. Then the tree $T_1$ induced
by the edges in $\{xx',x'y,x'x'',x''z\}$ and the trees $T_{i-1}$
induced by the edges in $\{xv_i,yv_i',zv_i'',v_iv_i',v_i'v_i''\} \
(2\leq i\leq s)$ are $s-1$ internally disjoint pendant $S$-Steiner
trees. Let $x_j,y_j$ be the vertices corresponding to $x,y$ in
$H(u_j) \ (4\leq j\leq r)$. The the trees $T_{j}'$ induced by the
edges in $\{xx_j, yy_j,zy_j,x_jy_j\} \ (4\leq j\leq r)$ are $r-3$
internally disjoint pendant $S$-Steiner trees. These trees together
with the trees $T_1,T_2,\cdots,T_{s-1}$ are $r+s-4$ internally
disjoint pendant $S$-Steiner trees.

Suppose that $x,y',z'$ are the same vertex in $H(u_1)$. Without loss
of generality, let $x=y'=z'$, $x=v_1$ in $H(u_1)$, $y=v_1'$ in
$H(u_2)$ and $z=v_1''$ in $H(u_3)$. Then the tree $T_{i-1}$ induced
by the edges in $\{xv_i,yv_i',zv_i'',v_iv_i',v_i'v_i''\} \ (2\leq
i\leq s)$ are $s-1$ internally disjoint pendant $S$-Steiner trees.
Let $x_j$ be the vertices corresponding to $x$ in $H(u_j) \ (4\leq
j\leq r)$. The the trees $T_{j}'$ induced by the edges in
$\{xx_j,yx_j,zx_j\} \ (4\leq j\leq r)$ are $r-3$ internally disjoint
pendant $S$-Steiner trees. These trees together with the trees
$T_1,T_2,\cdots,T_{s-1}$ are $r+s-4$ internally disjoint pendant
$S$-Steiner trees.

From the above argument, we conclude that for any
$S=\{x,y,z\}\subseteq V(G\Box H)$, there exist $r+s-3$ internally
disjoint pendant $S$-Steiner trees, and hence $\tau(S)\geq r+s-4$.
From the arbitraries of $S$, we have $\tau_k(G)=r+s-4$, as desired.
\qed
\end{pf}

From Lemma \ref{lem4-4}, if $G=K_{p+2,q-p+2}$, then
$\tau_3(L(G))=\tau_3(K_{p+2}\Box K_{q-p+2})=(p+2)+(q-p+2)-4=q$. So
the following result holds.
\begin{thm}\label{pro4-1}
For any two integers $p,q$ with $p\geq 2q$, there exists a graph $G$
such that $\tau_3(G)=p$ and $\tau_3(L(G))=q$.
\end{thm}


\begin{thebibliography}{1}

\bibitem{Bauer}
D. Bauer, R. Tindell, \emph{Graphs with prescribed connectivity and
line graph connectivity}, J. Graph Theory 3(1979), 393--395.

\bibitem{bondy}
J.A. Bondy, U.S.R. Murty, \emph{Graph Theory}, GTM 244, Springer,
2008.

\bibitem{DayA}
K. Day and A.-E. Al-Ayyoub, \emph{The cross product of
interconnection networks}, IEEE Trans. Parallel and Distributed
Systems 8(2)(1997), 109-118.


\bibitem{Capobianco}
M. Capobianco, J. Molluzzo, \emph{Examples and counterexamples in
Graph Theory}, North-Holland, Amsterdam 1978.

\bibitem{Chang}
S. Chang, \emph{The uniqueness and nonuniqueness of the triangular
association scheme}, Sci. Record 3(1959), 604--613.


\bibitem{Chartrand1}
G. Chartrand, S.F. Kappor, L. Lesniak, D.R. Lick, \emph{Generalized
connectivity in graphs}, Bull. Bombay Math. Colloq. 2(1984), 1--6.

\bibitem{Chartrand2}
G. Chartrand, F. Okamoto, P. Zhang, \emph{Rainbow trees in graphs
and generalized connectivity}, Networks 55 (4) (2010), 360--367.

\bibitem{Steeart} G. Chartrand, M. Steeart, \emph{The connectivity of
line graphs}, Math. Ann. 182(1969), 170--174.


\bibitem{Dirac}
G.A. Dirac, \emph{In abstrakten Graphen vorhandene vollst\"{a}ndige
4-Graphen und ihre Unterteilungen}, Math. Nach 22(1960), 61--85.

\bibitem{Du}
D. Du, X. Hu, \emph{Steiner tree problems in computer communication
networks}, World Scientific, 2008.


\bibitem{Fragopoulou}
P. Fragopoulou, S.G. Akl, \emph{Edge-disjoint spanning trees on the
star network with applications to fault tolerance}, IEEE Trans.
Computers 45(2)(1996), 174--185.

\bibitem{Grotschel1}
M. Gr\"{o}tschel, \emph{The Steiner tree packing problem in $VLSI$
design}, Math. Program. 78 (1997), 265--281.

\bibitem{Grotschel2}
M. Gr\"{o}tschel, A. Martin, R. Weismantel, \emph{Packing Steiner
trees: A cutting plane algorithm and commputational results}, Math.
Program. 72 (1996), 125--145.

\bibitem{Hager}
M. Hager, \emph{Pendant tree-connectivity}, J. Combin. Theory 38
(1985), 179--189.

\bibitem{Hedetniemi}
S.M. Hedetniemi, S.T. Hedetniemi, A.L. Liestman, \emph{A survey of
gossiping and broadcasting in communication networks}, Networks
18(1988), 1240--1268.


\bibitem{Hoffman} A. Hoffman, \emph{On the exeptional case in the
characterization of the arcs of a complete graph}, IBM J. Res. Dev.
4(1960), 487--496.

\bibitem{Jalote} P. Jalote, \emph{Fault Tolerance in Distributed Systems}, Prentice-Hall,
Englewood Cliffs, NJ, 1994.

\bibitem{Day}
K. Day, A.-E. Al-Ayyoub, \emph{The Cross Product of Interconnec-
tion Networks}, IEEE Trans. Parallel and Distributed Systems
8(2)(1997), 109--118.

\bibitem{Ku}
S. Ku, B. Wang, T. Hung, \emph{Constructing edge- disjoint spanning
trees in product networks}, Parallel and Distributed Systems, IEEE
Transactions on parallel and disjoited systems 14(3)(2003),
213--221.

\bibitem{LLM}
H. Li, X. Li, Y. Mao, \emph{On extremal graphs with at most two
internally disjoint Steiner trees connecting any three vertices},
Bull. Malays. Math. Sci. Soc. (2)37(3)(2014), 747--756.


\bibitem{LLSun}
H. Li, X. Li, Y. Sun, \emph{The generalied $3$-connectivity of
Cartesian product graphs}, Discrete Math. Theor. Comput. Sci. 14 (1)
(2012), 43--54.


\bibitem{LL}
S. Li, X. Li, \emph{Note on the hardness of generalized
connectivity}, J. Combin. Optimization 24 (2012), 389--396.

\bibitem{LLZ}
S. Li, X. Li, W. Zhou, \emph{Sharp bounds for the generalized
connectivity $\kappa_3(G)$}, Discrete Math. 310 (2010), 2147--2163.

\bibitem{LM1}
X. Li, Y. Mao, \emph{On extremal graphs with at most $\ell$
internally disjoint Steiner trees connecting any $n-1$ vertices},
Graphs \& Combin. 31(6)(2015), 2231--2259.

\bibitem{LM2}
X. Li, Y. Mao, \emph{The generalied $3$-connectivity of
lexigraphical product graphs}, Discrete Math. Theor. Comput. Sci.
16(1)(2014), 339--354.


\bibitem{LM3}
X. Li, Y. Mao, \emph{Nordhaus-Gaddum-type results for the
generalized edge-connectivity of graphs}, Discrete Appl. Math.
185(2015), 102--112.

\bibitem{LM4}
X. Li, Y. Mao, \emph{Graphs with large generalized
(edge-)connectivity}, accepted by Discuss. Math. Graph Theory,
arXiv: 1305.1089 [math.CO] 2013.


\bibitem{LMS}
X. Li, Y. Mao, Y. Sun, \emph{On the generalized (edge-)connectivity
of graphs}, Australasian J. Combin. 58(2)(2014), 304--319.

\bibitem{MaoL}
Y. Mao, H. Lai, \emph{On the pendant-tree connectivity of graphs},
submitted.

\bibitem{MaoLai}
Y. Mao, H. Lai, Graphs with given pendant-tree connectivity,
submitted.

\bibitem{Oellermann1}
O.R. Oellermann, \emph{Connectivity and edge-connectivity in graphs:
A survey}, Congessus Numerantium 116 (1996), 231--252.

\bibitem{Oellermann2}
O.R. Oellermann, \emph{On the $\ell$-connectivity of a graph},
Graphs \& Combin. 3(1987), 285--299.

\bibitem{Oellermann3}
O.R. Oellermann, \emph{A note on the $\ell$-connectivity function of
a graph}, Congessus Numerantium 60 (1987), 181--188.

\bibitem{Okamoto}
F. Okamoto, P. Zhang, \emph{The tree connectivity of regular
complete bipartite graphs}, J. Combin. Math. Combin. Comput. 74
(2010), 279--293.


\bibitem{Sherwani}
N.A. Sherwani, \emph{Algorithms for $VLSI$ Physical Design
Automation}, 3rd Edition, Kluwer Acad. Pub., London, 1999.


\bibitem{Ramanathan}
P. Ramanathan, D.D. Kandlur, K.G. Shin, \emph{Hardware-assisted
software clock synchronization for homogeneous distributed systems},
IEEE Trans. Comput. 39(1990), 514--524.

\bibitem{West} D. West, \emph{Introduction to Graph Theory (Second edition)},
Prentice Hall, 2001.

\bibitem{Whitney}
H. Whitney, \emph{Congruent graphs and the connectivity of graph},
Amer. J. Math. 54(1)(1932), 150--168.
\end{thebibliography}
\end{document}